\documentclass[english]{article}
\usepackage[T1]{fontenc}
\usepackage[latin9]{inputenc}
\usepackage{amsmath}
\usepackage{amssymb}
\usepackage{graphicx}
\usepackage[authoryear]{natbib}

\makeatletter

\providecommand{\tabularnewline}{\\}

\newcommand{\lyxaddress}[1]{
	\par {\raggedright #1
	\vspace{1.4em}
	\noindent\par}
}

\usepackage{amsmath}
\usepackage{accents}
\newlength{\dhatheight}

\makeatother

\usepackage{babel}
\begin{document}
\title{Generalized Estimators, Slope, Efficiency, and Fisher Information
Bounds}
\author{Paul Vos}
\maketitle
\begin{abstract}
Point estimators may not exist, need not be unique, and their distributions
are not parameter invariant. Generalized estimators provide distributions
that are parameter invariant, unique, and exist when point estimates
do not. Comparing point estimators using variance is less useful when
estimators are biased. A squared slope $\Lambda$ is defined that
can be used to compare both point and generalized estimators and is
unaffected by bias. Fisher information $I$ and variance are fundamentally
different quantities: the latter is defined at a distribution that
need not belong to a family, while the former cannot be defined without
a family of distributions, $M$. Fisher information and $\Lambda$
are similar quantities as both are defined on the tangent bundle $T\!M$
and $I$ provides an upper bound, $\Lambda\le I$, that holds for
all sample sizes -- asymptotics are not required. Comparing estimators
using $\Lambda$ rather than variance supports Fisher's claim that
$I$ provides a bound even in small samples. $\Lambda$-efficiency
is defined that extends the efficiency of unbiased estimators based
on variance. While defined by the slope, $\Lambda$-efficiency is
simply $\rho^{2}$, the square of the correlation between estimator
and score function.
\end{abstract}

\lyxaddress{\textbf{Key words:} nonasymptotics, generalized estimators, slope,
Fisher information, efficiency}

\section{Introduction}

Differential geometry has been used extensively in the study of higher
order asymptotic inference. While there are many references that could
be cited here, \citet{Amari1990-ly} stands out for its clear description
of the role that geometry plays in understanding and simplifying many
areas of statistical inference. This paper's finite sample results
are a departure from asymptotics, but we hope it continues geometry's
role in providing simplicity and intuition for inferential methods
and concepts. 

\section{Estimation}

\subsection{Two Views of Estimation}

One approach to estimation is that of choosing a point estimator with
good properties and reporting this estimate along with its standard
error. We label this the \emph{point-estimate-plus-se} approach. The
other approach is based on hypothesis testing where there is not just
one null hypothesis but a continuum. This is the \emph{continuum-of-hypotheses}
approach and it is the view expressed by \citet{Fisher1955} 
\begin{quote}
It may be added that in the theory of estimation we consider a continuum
of hypotheses each eligible as null hypothesis, and it is the aggregate
of frequencies calculated from each possibility in turn as true $\ldots$
which supply the likelihood function, $\ldots$ and other indications
of the amount of information available. 
\end{quote}
In the second approach, the point estimate, a real number, is replaced
with a function on the parameter space, a generalized estimate, that
is obtained from the continuum of hypotheses and estimation is described
in terms of properties of this function. 

\subsection{The Estimation Problem}

The estimation problem considered in this paper is described as follows.
There is a finite population having a distribution on values $\left\{ b_{1},b_{2},\ldots,b_{I}\right\} \subset\mathbb{R}$,
that, without loss of generality (wlog), we can order so that $b_{1}<b_{2}<\cdots<b_{I}$.
This population can be an actual population, such as all adults living
in North Carolina on January 1, 2022, or it can be hypothetical, such
as the collection of all pills that could be produced from a specific
production line. Whether actual or hypothetical, we call this the
real world population. The multivariate hypergeometric distribution
for a sample of size one along with the numeric labels is an exact
model for the data obtained from the real world population that we
denote as $\mathfrak{m}$. There is a one-dimensional family of distributions,
$M$, where each distribution in $M$ is dominated by a measure having
support $\mathcal{X}\subset\mathbb{R}$ such that $\left\{ b_{1},b_{2},\ldots,b_{I}\right\} \subset\mathcal{X}$.
The assumption is that there is at least one model $m\in M$ so that
$m$ is a suitable approximation to $\mathfrak{m}$. What constitutes
a suitable approximation will depend on the application. It may be
that the area histogram for $\mathfrak{m}$ is well approximated by
the area histogram of $m$ if $\mathcal{X}$ is countable or by the
graph of the density function if $\mathcal{X}$ is uncountable. \emph{}

With this assumption the estimation problem becomes that of using
points $y=(x_{1},x_{2},\ldots,x_{n})\in\mathcal{X}^{n}=\mathcal{Y}$
to distinguish models in $M$. To accommodate reduction to a sufficient
statistic and inference conditional on an ancillary statistic, the
space $\mathcal{Y}$ is adjusted accordingly. The real world population
and the exact model $\mathfrak{m}$ motivate the estimation problem
but neither play a role in defining or assessing estimators. Pointedly,
we will not take hypothetical repeated samples from either the real
world population or from $\mathfrak{m}$. We assume that the family
of distributions $M$ has the structure that makes $M$ a Riemannian
manifold and there is a global coordinate chart $\theta:M\rightarrow\Theta$
where $\Theta$ is an open interval in $\mathbb{R}$. 

\subsection{A Continuum of Hypotheses}

Fisher's continuum-of-hypotheses approach to estimation is less well-known
so we provide a few details. We are given a sample $y$ and a family
of models $M$ and the estimation problem is that of assessing, for
each model, the extent to which $y$ is an extreme observation. For
a fixed model, say the model with $\theta=\theta'$, the notion of
`extreme' is made precise by a test statistic which defines an ordering
$\preccurlyeq$ on $\mathcal{Y}$ given by 
\begin{equation}
y_{1}\preccurlyeq y_{2}\mbox{ iff }T(y_{1})\le T(y_{2}).\label{eq:ordering}
\end{equation}
This allows us to talk about a particular $y\in\mathcal{Y}$ as being
extreme, i.e., the value $T(y_{{\rm }})$ is extreme. Under mild regularity
conditions, $T^{-1}(t)$, for $t\in T(\mathcal{Y})$, provides a continuous
partition of $\mathcal{Y}$ into ordered subsets of co-dimension 1.
Since the ordering does not change if $T$ is replaced with $aT+b$
for $a,b\in\mathbb{R}$ and $a>0$ we can, wlog, define $T$ such
that $ET=0$ and $VT=1$.With $T$ standardized, $T(y_{{\rm }})$
is the number standard deviations $y_{{\rm }}$ is from 0, the mean
of the sampling distribution. 

Since we are considering a continuum of hypotheses we have a continuum
of test statistics $T(\cdot,\theta)$ where $T(y,\theta)$ indicates
how far $y$ is in the tails of the sampling distribution for each
model indexed by $\theta\in\Theta$. In comparing two tests $T_{1}$
and $T_{2}$, there would be a preference for $T_{1}$ if $T_{1}(y,\theta)$,
as a function of $\theta$, changed more rapidly than $T_{2}(y,\theta)$.
The derivative, which is the slope of the graph of $T(y,\theta)$,
measures this change at $\theta$. This discussion motivates the following
heuristic: `$T_{1}$ has more information than $T_{2}$ for distinguishing
models near $\theta$'. The slope of the test can be a function of
$y$ so we use the mean slope, $E_{\theta}T'(Y,\theta)$. Even though
the heuristic does not invoke a formal definition of `information',
it is worth noting that when $T$ is the score, the mean slope is
the Fisher information and when $T$ is the standardized score the
mean slope is the square root of the Fisher information. 

The notion of a continuum of tests is formalized below using a generalized
estimator and this approach to estimation can be summarized in the
following desideratum: 
\begin{quote}
Generalized estimators should have large\footnote{Here, 'large' means in absolute value.}
mean slope. 
\end{quote}
This is related to, but distinct from, the desideratum for point estimation:
\begin{quote}
Estimators should have small standard deviation. 
\end{quote}
Estimators are often compared in terms of variance, rather than standard
deviation, and squared mean slope is used to show the connection between
these two approaches to estimation. However, the geometric intuition
for understanding the relationship between estimators, both point
and generalized estimators, is best described using the mean slope
rather than its square. 

\section{Generalized Estimators}

Estimators, both point and generalized, can be described by their
aggregate or distributional properties on $\mathcal{Y}$. Such properties
can be used to compare and select estimators before any specific data
have been obtained. We begin by considering aggregate properties,
that is \emph{estimators}. Section \ref{sec:Confidence-Intervals-as}
describes the inferential role of the generalized \emph{estimate}
obtained from a single sample $y\in\mathcal{Y}$.

\subsection{Space of Generalized Estimators $\mathcal{G}$ \label{subsec:Space-of-Generalized}}

A point estimator for $\theta$ is a function $u:\mathcal{Y}\rightarrow\Theta$.
A generalized estimator for $\theta$ is a function $g:\mathcal{Y}\rightarrow C^{k}\left(\Theta\right)$,
the space of functions on $\Theta$ whose $k^{th}$ derivative is
continuous. We let $C(\Theta)$ be the space of functions of required
smoothness. For fixed $y$, point estimation $y\mapsto u(y)=\hat{\theta}\in\Theta$
while generalized estimation $y\mapsto g_{y}(\theta)\in C\left(\Theta\right)$.
When $y$ is fixed we simplify the notation by dropping the subscript,
$g(\theta)=g_{y}(\theta)$. The context will indicate whether $g=g(\theta)$
is an estimate (fixed $y$) or an estimator (distribution of $y$
values). This follows the notational simplification for point estimation
where $\hat{\theta}$ may refer to an estimate or to an estimator
and expressions such as $E\hat{\theta}$ indicate that $\hat{\theta}$
is an estimator. Likewise, $Eg$ indicates that $g$ is an estimator
and if the expectation were defined using an integral we may use the
full notation $g_{y}(\theta)$ or even $g(\theta;y)$. The generalized
estimate needs to be a function on $\Theta$ for it to represent the
continuum-of-hypotheses approach since each hypothesis specifies a
point in $\Theta$. For the null hypothesis $H_{\circ}:\theta=\theta_{\circ}$,
$g(\theta_{\circ})$ is the value of a test statistic for the sample
$y$.   The function $g(\theta)$ shows how the test statistic for
fixed $y$ changes over the continuum of hypotheses indexed by $\theta$. 

For another parameterization $\xi$, $y\mapsto g(\xi)\in C\left(\Xi\right)$
and while $g(\theta)$ and $g(\xi)$ are different functions, they
take the same value when $\theta$ and $\xi$ label the same point
$m\in M$. Using more precise notation, a generalized estimate for
$y$ is a function $g_{M}:M\rightarrow\mathbb{R}$ that when expressed
in the $\theta$ parameter is $g_{\Theta}:\Theta\rightarrow\mathbb{R}$
and $g_{\Xi}:\Xi\rightarrow\mathbb{R}$ in the $\xi$ parameter. For
any point $m_{\circ}\in M$, $g_{M}(m_{\circ})$, $g_{\Theta}(\theta(m_{\circ}))$,
and $g_{\Xi}(\xi(m_{\circ}))$ all take the same value. Furthermore,
since $y$ was any point in $\mathcal{Y}$ and $m_{\circ}$ specifies
a distribution on $\mathcal{Y}$, the estimators $g_{M}(m_{\circ})$,
$g_{\Theta}(\theta_{\circ})$, and $g_{\Xi}(\xi_{\circ})$ all have
the same distribution. It is this property that make generalized estimators
parameter invariant. 

For point estimation the function $u$ must be measurable and we generally
require $E_{\theta}u^{2}<\infty$ for all $\theta$. For the generalized
estimator we require that $g(\theta)$ is measurable for all $\theta$.
The moment assumptions are that $E_{\theta}g(\theta)=0$, $V_{\theta}g(\theta)<\infty$,
and $E_{\theta}(g(\theta)\ell'(\theta))\ge0$ where $\ell(\theta)$
is the log likelihood for $M$. The smoothness assumptions are that
$\partial_{\theta}V_{\theta}g(\theta)$ is finite and $g(\theta)$
can be differentiated under the integral so that $\partial_{\theta}E_{\theta}g(\theta)=0$
yields the identity
\begin{equation}
-Eg'(\theta)=E\left(g(\theta)\ell'(\theta)\right)\label{eq:FundamentalEQ1}
\end{equation}
where the subscript $\theta$ on $E$ has been dropped since the expectations
will always be taken with respect to the same value for the parameter
as appears in the functions. Further notational simplification emphasizes
that we are interested in a continuum of values
\begin{equation}
-Eg'=E(g\ell')\label{eq:FundamentalEQ1simple}
\end{equation}
where for a.e. $y$, $g',\ g,\ \ell'\in C(\Theta)$ and the expectations
also are in $C(\Theta)$.  Any function $g$ satisfying these conditions
is a generalized estimator and we denote the space of all generalized
estimators as $\mathcal{G}(M,\mathcal{Y})$ which we abbreviate to
$\mathcal{G}$ when the context is clear regarding the model space
and sampling distribution. Since $-g(\theta)$ simply reverses the
ordering of $\mathcal{Y}$ specified by $g(\theta)$ we don't want
both of these functions in $\mathcal{G}$; the restriction $E(g\ell')\ge0$
means $\ell'\in\mathcal{G}$. 

The relation $g_{1}\sim g_{2}$ if-f $g_{1}(\theta)=k(\theta)g_{2}(\theta)$
for some $k\in C(\Theta)$ with $k(\theta)>0$ is an equivalence relation.
The equivalence of $g_{1}$ and $g_{2}$ is motivated by the fact
that for each $\theta_{\circ}$, $g_{1}(\theta_{\circ})$ and $g_{2}(\theta_{\circ})$
provide the same ordering of $\mathcal{Y}$ so that the tail area
obtained for each $y\in\mathcal{Y}$ is the same for both $g_{1}$
and $g_{2}$. 

Inferential properties of equivalence class $[g]$ are conveniently
represented by the standardized estimator 
\[
\bar{g}=g/\sqrt{V(g)}.
\]
The standardized score $\ell'$ will be denoted as
\[
\hat{s}=\ell'/\sqrt{V(\ell')}\in[\ell'].
\]

A generalized estimator $g(\theta)$ defines a point estimator $\hat{\theta}_{g}$
by the equation $g(\hat{\theta}_{g})=0$. Note that $\hat{\theta}_{g}$
need not be unique and it need not exist.

\subsection{Mapping Point Estimators to $\mathcal{H}\subset\mathcal{G}$}

If $u:\mathcal{Y}\rightarrow\Theta$ is a point estimator and $\upsilon(\theta)=Eu$,
then
\begin{equation}
h(\theta)=u(y)-\upsilon(\theta)\label{eq:h from u}
\end{equation}
is a generalized estimator. If $\hat{\theta}=u(y)$ is unbiased for
$\theta$, then $\hat{\theta}_{h}=\hat{\theta}$. We use the notation
$h$ for generalized estimators that are the difference between two
functions, one having domain $\mathcal{Y}$ and the other $\Theta$.
All graphs for $h(\theta)$ have the same shape being the function
$-\upsilon(\theta)$ shifted by the amount $u(y)$. For $\hat{\theta}=u(y)$
unbiased $\upsilon(\theta)=\theta$ so the graphs of $h(\theta)$
are straight lines with slope $-1$; in general, the graphs of $\bar{h}(\theta)$
will not be straight lines but the shapes will be the same. The set
of generalized estimators of this form is $\mathcal{H}\subset\mathcal{G}.$
Note that $h$ is parameter-invariant since 
\[
h(\xi)=u(y)-\upsilon\circ\theta_{\xi}(\xi)
\]
where $\theta_{\xi}:\Xi\rightarrow\Theta$. Also, for $k_{1},k_{2}\in\mathbb{R}$
with $k_{1}>0$, $k_{1}u+k_{2}\in[u]=[u-Eu]$. Equivalence of statistic
$u$ under affine transformations is required for inference to be
invariant with respect to units on $\mathcal{Y}$ (e.g., degrees Fahrenheit
or Centigrade). The group of transformation defining equivalence for
generalized estimator $g$ provides invariance with respect to parameterizations
on $M$. 

The point estimator $\hat{\theta}=u(y)$ is not parameter-invariant
since its distribution need not be the same as that of $\hat{\xi}=\theta_{\xi}^{-1}(\hat{\theta})$.
Point \emph{estimators} for different parameterizations produce different
generalized estimators even when the point \emph{estimates} are parameter-invariant.
An example appears in Section \ref{subsec:Maximum-Likelihood-Estimator}.

\subsection{Squared Slope $\Lambda(g)$}

The squared slope of $g\in\mathcal{G}$ is a function $\Lambda(g)\in C(\Theta)$
defined by
\begin{align*}
\Lambda(g) & =\left(E\bar{g}'\right)^{2}\\
 & =(Eg')^{2}/V(g).
\end{align*}
Clearly, $\Lambda(g_{1})=\Lambda(g)$ whenever $g_{1}\sim g$. The
slope for a point estimator is defined as the slope of its corresponding
generalized estimator. The squared slope for $\hat{\theta}=u(y)$
is that of $h=\hat{\theta}-\upsilon(\theta)$
\begin{align*}
\Lambda(\hat{\theta}) & =(\upsilon')^{2}/V(\hat{\theta}).
\end{align*}
If $\hat{\theta}$ is unbiased, then $\upsilon(\theta)=\theta$ and
\begin{equation}
\Lambda(\hat{\theta})=1/V(\hat{\theta}).\label{eq:slopeUnbiasedPtEstimate}
\end{equation}

A caution regarding terminology: the slope of a generalized estimator
is a mean since there is a distribution of curves of varying slopes;
the slope of $\hat{\theta}-\upsilon(\theta)$ is $-\upsilon(\theta)'$
since each curve has the same slope. 

\section{Fisher Information Bounds}

Fisher information provides an upper bound for squared slope $\Lambda$
and a lower bound for the variance, $V$, of an estimator. In the
next section the bounds are calculated and efficiencies defined by
$\Lambda$ and $V$ are compared. Section \ref{subsec:Var_is_Pointwise}
shows that $V$ and Fisher information $I$ are fundamentally different;
variance of a point estimator is not geometric while $\Lambda$ and
$I$ are tensors. 

\subsection{Bounds for $\Lambda$ and $V$}

The Fisher information for $M$ is $I_{1}=-E\ell_{1}''(\theta,x)$
where $\ell_{1}$ is the log likelihood for models in $M$; $I_{1}$
represents the average amount of information in a single observation
$x$. The Fisher information in $y$ is $I=I_{\mathcal{Y}}=-E\ell''(\theta,y)$
where $\ell$ is the log likelihood for $y$. For $\mathcal{Y}=\mathcal{X}^{n}$,
$\ell=\sum\ell_{1}(\theta,x_{i})$ and $I=nI_{1}$. The Fisher information
upper bound for $\Lambda$ is obtained from $-Eg'=E(g\ell')$. Expressing
the covariance $E(g\ell')$ in terms of the correlation yields
\begin{equation}
-Eg'=\rho(g,\ell')\sqrt{V(g)V(\ell')}.\label{eq:CorrEqn}
\end{equation}
Squaring this last equation and making substitutions $\Lambda(g)=(Eg')^{2}/V(g)$
and $I=V(\ell')$ shows
\begin{equation}
\Lambda(g)=\rho^{2}(g,\ell')I\label{eq:FishSlopeBound}
\end{equation}
so that the Fisher information $I$ is the upper bound for $\Lambda(g)$
since $\rho\le1$. Clearly, $\Lambda(\ell')=I$ so that the score
achieves the squared slope upper bound. Since $\Lambda(\hat{\theta})=\Lambda(\hat{\theta}-\upsilon(\theta))$,
equation (\ref{eq:FishSlopeBound}) provides a bound for point estimator
$\hat{\theta}$, 
\begin{equation}
\Lambda(\hat{\theta})=\rho^{2}(\hat{\theta},\ell')I.\label{eq:FishSlopeBound2}
\end{equation}

For the variance bound, we substitute $\hat{\theta}-\upsilon(\theta)$
for $g$ in (\ref{eq:CorrEqn}) and solving for variance yields 
\begin{equation}
V(\hat{\theta})=\frac{(\upsilon')^{2}}{\rho^{2}(\hat{\theta},\ell')I}\ge(\upsilon')^{2}I^{-1}\label{eq:FishVarBound}
\end{equation}
where the equality in (\ref{eq:FishVarBound}) holds provided $\rho\not=0$.
The bound on $\Lambda(\hat{\theta})$ has no restriction on $\rho$.
When $\hat{\theta}$ is unbiased $\upsilon(\theta)=\theta$ so that
(\ref{eq:FishVarBound}) gives the Cram\'{e}r-Rao lower bound $V(\hat{\theta})\ge I^{-1}$.
The bound on $\Lambda(\hat{\theta})$ is the same whether $\hat{\theta}$
is biased or unbiased. 

The \emph{$V$-efficiency} of an unbiased point estimator $\hat{\theta}$
is defined in terms of the ratio $V(\hat{\theta})$ and the variance
lower bound
\begin{equation}
\mbox{Eff}^{V}(\hat{\theta})=I^{-1}/V(\hat{\theta)}.\label{eq:Veff}
\end{equation}
We define the \emph{$\Lambda$-efficiency} of any generalized estimator
$g$ using the squared slope and its upper bound 
\begin{equation}
\mbox{Eff}^{\Lambda}(g)=\Lambda(g)/I=\rho^{2}(g,\ell').\label{eq:LambdaEff}
\end{equation}
Setting $g=\hat{\theta}-\upsilon(\theta)$, the $\Lambda$-efficiency
for a point estimate is, like that of generalized estimator $g$,
simply the square of the correlation between the estimator and the
score
\[
\mbox{Eff}^{\Lambda}(\hat{\theta})=\rho^{2}(\hat{\theta},\ell').
\]
For unbiased $\hat{\theta}$ and $\rho\not=0$, the equality in (\ref{eq:FishVarBound})
shows $\mbox{Eff}^{V}(\hat{\theta})=\mbox{Eff}^{\Lambda}(\hat{\theta})$
so that $\Lambda$-efficiency is an extension of $V$-efficieny to
unbiased point estimators. In fact, $\Lambda$-efficiency extends
$V$-efficiency to all generalized estimators.

Since $I=nI_{1}$ when $\mathcal{Y}$ is $\mathcal{X}^{n}$, $\Lambda(g)=\rho^{2}nI_{1}$
and $\Lambda$-efficiency can be interpreted in terms of an effective
sample size
\[
n(g)=\rho^{2}n
\]
where $\rho$ is the correlation between $g$ and $\ell'$. Or, the
$\Lambda$-efficiency loss of using $g$ rather than $\ell'$ is a
loss of $(1-\rho^{2})n$ observations.

\subsection{Variance is Defined Point-wise \label{subsec:Var_is_Pointwise}}

For any statistic $u:\mathcal{Y}\rightarrow\mathbb{R}$ , the variance
of $u$ for model $m$ depends only on $m$; it does not depend on
the family $M$ that contains $m$. Fisher information is defined
on the tangent bundle $T\!M$ and so cannot be defined for an isolated
point outside of $M$. Fisher information and variance are fundamentally
different quantities. Like Fisher information, $\Lambda$ is defined
on $T\!M$ and the Fisher information bound, $\Lambda=\rho^{2}I$,
provides interpretations for finite samples. We illustrate the difference
between the variance and squared slope $\Lambda$ of an estimator
using an example where a model $m$ is not isolated but belongs to
two different manifolds. 

Let $\mathcal{X}=\left\{ \left(x_{1},x_{2}\right)\in\mathbb{R}^{2}\right\} $,
$M$ be the two dimensional manifold of bivariate normal distributions
with zero correlation and variance one, and $\mu$ be the mean parameterization.
Reduction of the sample space to the minimal sufficient statistic
provides the two dimensional sample space $\mathcal{Y}=\mathbb{R}^{2}$.
We consider submanifolds $M_{1}$ and $M_{2}$ that intersect in a
point $m_{\circ}$. We start with $M_{1}$ consisting of all points
having $\mu=(\theta,0)$ and $M_{2}$ all points with $\mu=(0,\theta)$
so that $m_{\circ}$ is the point with $\mu=(0,0).$ In the expectation
parameter space, $M_{1}$ is the horizontal axis and $M_{2}$ is the
vertical axis. We compare two estimators: $\bar{x}_{j}=n^{-1}\sum x_{ji}$,
for $j\in\left\{ 1,2\right\} $. Both $\bar{x}_{1}$ and $\bar{x}_{2}$
have variance $1/n$ that equals the variance lower bound for all
models in $M_{1}$ and $M_{2}$. Both estimators are unbiased for
$\theta$ at $m_{\circ}$. In terms of bias and variance the estimators
are indistinguishable at $m_{\circ}$. For inference in the family
$M_{1}$, $\bar{x}_{1}$ has $\rho(\bar{x}_{1},\ell')=1$ while $\bar{x}_{2}$
has $\rho(\bar{x}_{2},\ell')=0$. The variance of an estimator at
a point does not distinguish between an optimal estimator and one
that is worthless (its distribution is the same for all models in
the manifold). 

In contrast, $\Lambda_{M_{1}}(\bar{x}_{1})=n$ attaining the slope
upper bound while $\Lambda_{M_{1}}(\bar{x}_{2})=0$. The squared slope
correctly distinguishes these estimators in terms of their inferential
properties. The difference is that $\Lambda_{M_{1}}(u)=\rho^{2}(u,\ell')n$
is a function of $\ell'\in T\!M_{1}$. 

This problem disappears when sufficiency arguments are used to replace
$\mathcal{Y}$ with $\mathcal{Y}_{j}=\left\{ y:y=n^{-1}\sum x_{ji}\right\} $
for $M_{j}$. However, by replacing $M_{1}$ and $M_{2}$ with manifolds
whose expectation parameters are circles of large radii intersecting
at $(0,0)$ with centers on the vertical and horizontal axes the minimal
sufficient statistic is two dimensional and the sample distributions
for models in $M_{1}$ and $M_{2}$ are defined on $\mathcal{Y}$.
The expectation and variance of $\bar{x}_{1}$ and $\bar{x}_{2}$
at $\mu=(0,0)$ are unaffected by changing the manifolds; $\Lambda_{M_{1}}(\bar{x}_{j})$
will be a function of $\theta$ that for $\theta=0$ approaches the
upper bound $n$ when $j=1$ and zero when $j=2$. 

The problem with using variance to assess an estimator is lessened
when the variance is considered a function on the manifold. The argument
that estimators should be assessed at all distributions is cogent,
if not obvious, but begs the question of why estimators should not
be allowed to depend on the distribution. 

\section{$\mathcal{G}$ and Confidence Intervals \label{sec:Confidence-Intervals-as}}

After selecting a generalized estimator $g$ based on its aggregate
properties, we use the resulting estimate $g(\theta)\in C(\Theta)$
corresponding to a particular value $y\in\mathcal{Y}$. We allow $g(\theta)$
to be any generalized estimate although $\Lambda$-optimality would
suggest using $\ell'(\theta)$. 

For $g(\theta)$ monotone on $\Theta$, the interval
\[
(\theta_{lo},\theta_{hi})=\left\{ \theta:-k<\bar{g}(\theta)<k\right\} 
\]
is the $g$-\emph{confidence interval for $\theta$ at level} $k$.
With $g=\ell'$ and $k=2$, $(\theta_{lo},\theta_{hi})=\hat{s}^{-1}(\pm2)$
consists of all models $m$ such that $y$ is within 2 standard deviations
of the most likely sample value when $\mathcal{Y}$ is ordered by
the score. The justification for ordering $\mathcal{Y}$ by $\ell'$
is that it has the greatest slope among all estimators in $\mathcal{G}$.
That is, $\ell'$ is $\Lambda$-optimal. This interpretation for $\hat{s}^{-1}(\pm2)$
holds for any sample size $n$; the interpretation describing the
percentile that $y$ is in $\mathcal{Y}$ will depend on the extent
to which the distribution of $\hat{s}(\theta)$ is approximated by
the standard normal for all $\theta\in(\theta_{lo},\theta_{hi})$. 

A linear approximation to the estimate $\bar{g}(\theta)$ is of the
form 
\[
b(\theta-\hat{\theta})
\]
where $b$ may depend on $y$ but is constant in $\theta$ and $\hat{\theta}$
satisfies $\bar{g}(\hat{\theta})=0$. Function $b(\theta-\hat{\theta})$
is considered a good approximation to $\bar{g}(\theta)$ when confidence
intervals obtained from the linear function are a suitable approximation
to those obtained from $\bar{g}(\theta)$ over a subset $\Theta_{1}\subset\Theta$.
That is, for $k\in\mathbb{R}$, 

\begin{align}
\left\{ \theta:-k<\bar{g}(\theta)<k\right\}  & \doteq\left\{ \theta:-k<b(\theta-\hat{\theta})<k\right\} \nonumber \\
 & =\left\{ \theta:\hat{\theta}+k/b<\theta<\hat{\theta}-k/b\right\} \nonumber \\
 & =\hat{\theta}\pm k|b|^{-1}.\label{eq:slopeCI}
\end{align}
Comparing (\ref{eq:slopeCI}) to the estimate-plus-standard-error
approach, $\hat{\theta}\pm k\widehat{se}$, shows that, for unbiased
$\hat{\theta}$, the standard error $\widehat{se}=\sqrt{V(\hat{\theta})}$
is the reciprocal of the absolute value of the slope $b$ of the linear
approximation to $\bar{g}(\theta)$. The linear approximation will
not be parameter invariant.

\section{A Related Class of Estimators $\mathcal{G}_{{\rm lrt}}$}

We consider the class of all functions $G:\mathcal{Y}\rightarrow C(\Theta)$
such that $G'\in\mathcal{G}$ and $K=\sup_{\theta\in\Theta}G$ is
finite and measurable. We denote the set of functions $G-K$ as $\mathcal{G}_{{\rm lrt}}$
since the log likelihood ratio statistic in this class plays the role
of the score in the class of estimators $\mathcal{G}$. Following
the same logic, the space $\mathcal{G}$ defined in Section \ref{subsec:Space-of-Generalized}
would be denoted $\mathcal{G}_{{\rm score}}$. Note that $\sup_{\theta}G=0$
for all $G\in\mathcal{G}_{{\rm lrt}}$. When $G_{\Theta}=G$ is expressed
in parameter $\xi$, $\partial_{\xi}G_{\Xi}=\left(\partial_{\xi}\theta_{\Xi}\right)\partial_{\theta}G_{\Theta}$
so that $G'_{\Xi}\sim G_{\Theta}'$. Consequently, the functions in
$\mathcal{G}_{{\rm lrt}}$, like those in $\mathcal{G}_{{\rm score}}$,
are parameter-invariant and so each parameterized version is obtained
from a function $G_{M}$ having domain $M$. Usually the context will
indicate the domain of the function so we drop the subscript notation
on $G$. 

The squared slope can be extended to $\mathcal{G}_{{\rm lrt}}$ by
defining $\Lambda(G)=\Lambda(G')$ where the interpretation is now
in terms of the Laplacian of $G$ rather than the slope of $G'$ .
Since $\hat{s}$ is $\Lambda$-optimal so is 
\begin{equation}
\widehat{S}(\theta)=2\left(\ell(\theta,y)-\sup_{\theta\in\Theta}\ell(\theta,y)\right)\label{eq:lrtShat}
\end{equation}
where the constant 2 is chosen for convenience as the distribution
of $-\widehat{S}(\theta)$ is approximately Chi-squared with 1 degree
of freedom. 

The choice between using $\hat{s}$ or $\widehat{S}$ must depend
on considerations other than $\Lambda$-optimality. Advantages of
$\hat{s}$ are that the ordering of $\mathcal{Y}$ accounts for the
direction of departures from the model and confidence intervals can
be interpreted in terms of standard deviations. When each model in
$M$ specifies a symmetric distribution on $\mathcal{Y}$, the tails
specified by $\hat{s}$ are identical (mirror images) so that $\widehat{S}$,
or any $G$, can be used while still accounting for the direction
of the departure. When $M$ is the family of normal location models,
confidence intervals from $\hat{s}$ and $\widehat{S}$ are identical.
An advantage of $\widehat{S}$ is that it can be used for inference
for two or more parameters considered jointly. 

\section{Examples}

\subsection{Bernoulli Family}

Consider a dichotomous attribute for individuals in a population.
We use the sample space $\mathcal{X}=\left\{ 0,1\right\} $ where
$1$ indicates the attribute is present and the statistical manifold
is
\[
M=\left\{ m:\mathcal{X}\rightarrow\mathbb{R}:0<m(1)<1,\ m(0)+m(1)=1\right\} .
\]
The sample space $\mathcal{X}^{n}$ can be reduced to $\mathcal{Y}=\left\{ 0,1,\ldots,n\right\} $
using the sufficient statistic $y=\sum x_{i}$. The log likelihood
function for parameter $p(m)=m(1)$ is 
\[
\ell_{P}(p)=y\log p+(n-y)\log(1-p)+k(n,y)
\]
with score
\[
\hat{s}_{P}(p)=n^{-1/2}\frac{y-np}{\sqrt{p(1-p)}}.
\]
For the log-odds parameterization $\theta=\theta(m)=\log(m(1)/m(0))$,
the score is
\[
\hat{s}_{\Theta}(\theta)=n^{-1/2}\frac{y-np_{\Theta}(\theta)}{\sqrt{p_{\Theta}(\theta)(1-p_{\Theta}(\theta))}}
\]
where $p_{\Theta}(\theta)=e^{\theta}/(e^{\theta}+1$).

\begin{figure}
\includegraphics[scale=0.9]{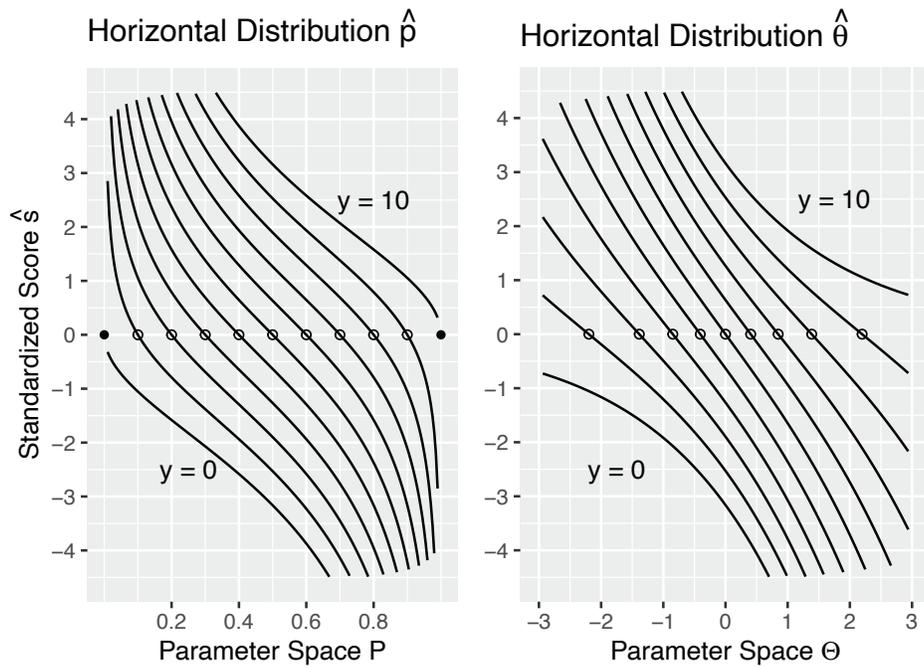}\caption{Maximum likelihood (ml) point estimate for parameters $p$ and $\theta$.
The nine open circles are the ml point estimates; the models in $M$
corresponding to the estimates in the $P$ are the same as those in
$\Theta$. Maximum likelihood provides parameter invariant \emph{estimates}.
\label{fig:BinomPoint}}

\end{figure}

\begin{figure}

\includegraphics[scale=0.9]{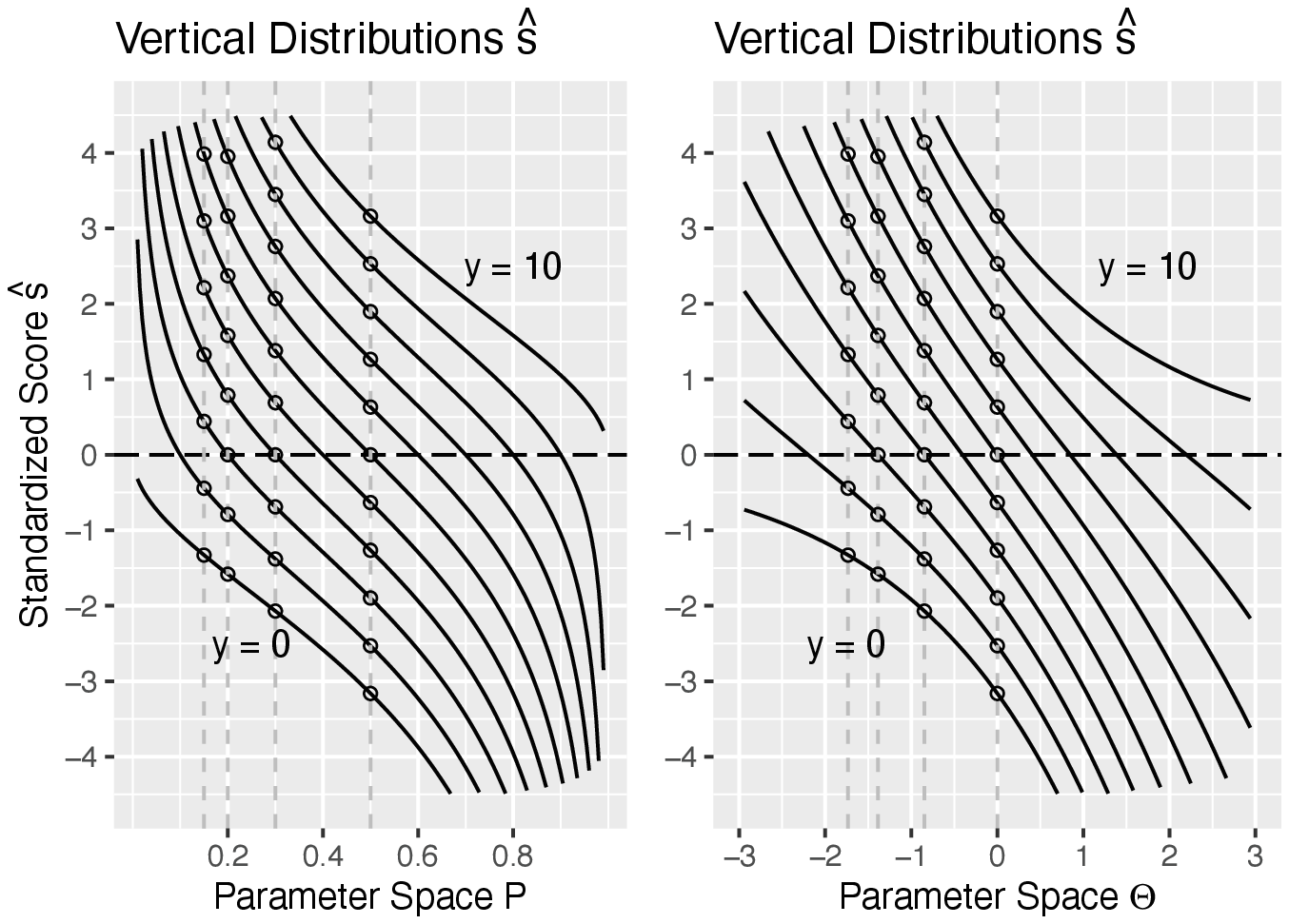}\caption{Generalized estimator $\hat{s}$. The two sets of eleven curves are
the generalized estimator $\hat{s}$ expressed in the $P$ and $\Theta$
parameter spaces. The distribution of $\hat{s}$ is shown for four
models in $M$ as vertical slices of the curves. The distribution
of $\hat{s}(p)|_{p=1/2}$ is identical to that of $\hat{s}(\theta)|_{\theta=0}$.
This is true for all $p$ and $\theta$ when both label the same model
in $M$. Generalized estimators provide parameter invariant \emph{estimators}.
\label{fig:BinomGen}}

\end{figure}

\subsubsection{Maximum Likelihood Estimator\label{subsec:Maximum-Likelihood-Estimator}}

For $y\not\in\left\{ 0,n\right\} $ the maximum likelihood (ml) point
estimates $\hat{p}=y/n$ and $\hat{\theta}=\theta_{P}(\hat{p})$ identify
the same point in $M$, namely, the distribution for which $m(1)=y/n$.
For $y\in\left\{ 0,n\right\} $ the ml point estimators are not defined
so bias and variance cannot be used to compare these estimators. We
extend the definition of maximum likelihood estimate in the next section
where the variance is finite for the parameter $p$ but infinite for
$\theta$. Figure \ref{fig:BinomPoint} displays the eleven standardized
score functions for both parameters when $\mathcal{Y}=\left\{ 0,1,\ldots,10\right\} $.
The location of the circles along the horizontal axis are the nine
ml point estimates that exist. 

For all $y\in\mathcal{Y}$ the generalized ml estimates $\hat{s}(p)$
and $\hat{s}(\theta)$ exist and for every $m_{\circ}\in M$ $\hat{s}(p_{\circ})=\hat{s}(\theta_{\circ})$
where $p_{\circ}=p(m_{\circ})$ and $\theta_{\circ}=\theta(m_{\circ})$
so the distribution of these estimators are parameter invariant. The
curves in the left panel of Figure \ref{fig:BinomGen} are the eleven
generalized estimates expressed using parameter $p$; the right panel
shows these generalized estimates using parameter $\theta$. These
are the same curves that appear in Figure \ref{fig:BinomPoint} but
we are now interested in the entire curve, not just where nine of
these cross the parameter axis. 

While the distribution for point estimators is restricted to the intersection
of these curves with the horizontal axis, generalized estimators describe
a continuum-of-hypotheses where each hypothesis considered in turn
is represented by where these curves intersect the vertical line specified
by a null value $p_{\circ}$ or $\theta_{\circ}.$ Figure \ref{fig:BinomGen}
displays the distribution for four hypotheses. Notice in each case
the corresponding vertical distributions are the same for both parameterizations. 

For any parameterization, these vertical distributions have the same
mean and same standard deviation for all parameter values, and so
provide little information for estimation. Estimation using this continuum
of vertical distributions is described by the mean slope. For example,
the mean of the eleven slopes at $\theta=0$ has absolute value $\sqrt{\Lambda(\hat{s}(\theta))}|_{\theta=0}$.

The Fisher information bound expressed in terms of the desideratum
that generalized estimators should have large mean slopes is described
as follows. For each $g\in\mathcal{G}$ there is a family of (at most)
eleven curves that are the graphs over a parameter space, say $\Theta$.
We seek an estimator $g$ such that the mean slope is large and the
Fisher information bound says its absolute value is at most $\sqrt{I(\theta)}$
for all $\theta\in\Theta$. Since mean slope and Fisher information
are tensors, Fisher information provides a bound for all parameterizations.

\subsubsection{Uniformly Minimum Variance Unbiased Estimators}

The estimator $\tilde{p}=y/n$ for $y\in\mathcal{Y}$ equals the ml
estimator for $y\not\in\left\{ 0,n\right\} $ and so is called the
\emph{extended} ml estimator. The extended ml estimator is the unique
uniformly minimum variance unbiased estimator (UMVUE) for parameter
$p$.

The parameter $\xi=p^{2}$ taking values in $\Xi=(0,1)$ has UMVUE
estimator $\tilde{\xi}=y(y-1)/(n(n-1))$ \citep[page 151]{Barnett1999-pf}.
There are two UMVU estimators, $\tilde{p}$ and $\tilde{\xi}$, and
since they are not both unbiased in the same parameter, we cannot
compare them directly in terms of their variances. We consider a third
estimator $(y/n)^{2}$ that is unbiased for $\eta=E(y/n)^{2}$. The
parameter $\eta$ is a function of the sample size so its estimator
is not likely to be selected based on its unbiasedness. 

\begin{figure}
\includegraphics[scale=0.9]{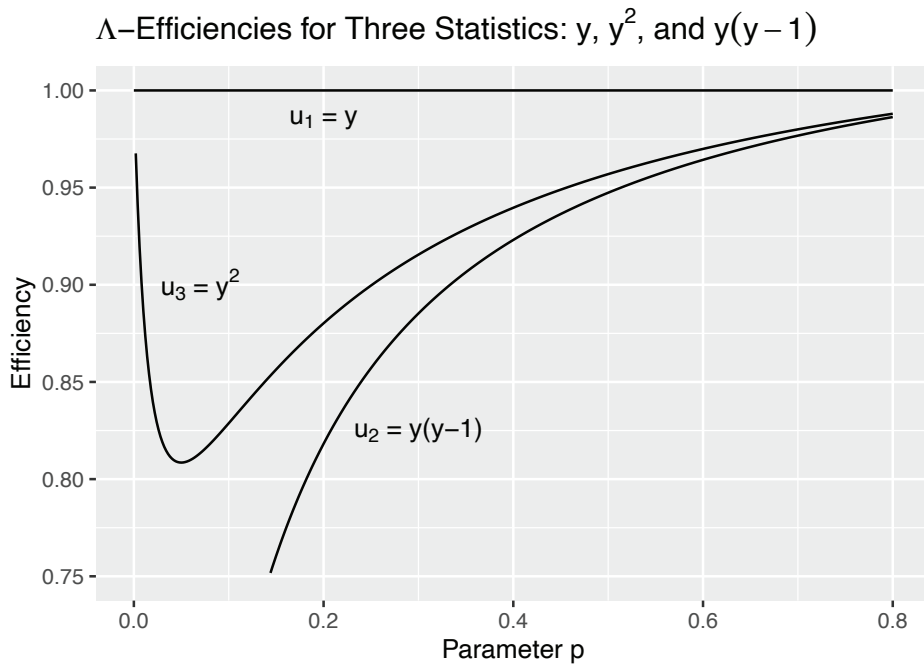}\caption{$\Lambda$-Efficiencies of three statistics for Bernoulli manifold
$M$ and sample space $\mathcal{Y}=\left\{ 0,1,\ldots,10\right\} $.
\label{fig:Efficiencies-of-three}}
\end{figure}

The generalized estimators associated with these point estimators
\[
h_{j}(\theta)=u_{j}-Eu_{j}
\]
are defined by the statistics $u_{1}=y$, $u_{2}=y(y-1)$, and $u_{3}=y^{2}$,
respectively. Note $\bar{h}_{1}=\hat{s}$. Figure \ref{fig:Efficiencies-of-three}
shows the $\Lambda$-efficiencies for these statistics. The parameter
$\theta$ can be any parameterization, parameter $p$ is used in Figure
\ref{fig:Efficiencies-of-three}. The curves will look different for
other parameterizations but the vertical slices giving the efficiencies
for corresponding points in the manifold are the same; that is, $\Lambda$-efficiency
can be defined directly on $M$. 

In terms of $\Lambda$-efficiency, the estimate for $p$ ($u_{1}=y)$
is fully efficient, the estimate for $\xi$ ($u_{2}=y(y-1)$) has
the lowest efficiency approaching zero for small $p$, and the estimate
based on $u_{3}=y^{2}$ has efficiency between the other statistics.
The fact that $u_{2}=y(y-1)$ is zero for both $y=0$ and $y=1$ is
consistent with efficiency approaching zero for models with $p$ near
zero. When $p<0.17$ using $\tilde{\xi}$ rather than $\tilde{p}$
results in an efficiency loss of at least 2 observations. Note that
$u_{1}$ and $u_{3}$ are both sufficient statistics but they are
not both $\Lambda$-efficient. 

Figure \ref{fig:Efficiencies-of-three} shows an asymmetry in the
estimator $\tilde{\xi}$ whose efficiency is close to one for $p$
near one but close to zero for $p$ near zero. For a dichotomous population
where one of the categories is rare, the efficiency of $\tilde{\xi}$
will depend on which category is labeled 1. In other words, statistics
$y(y-1)$ and $(n-y)(n-y-1)$ have very different efficiencies, even
though they describe the same inference problem with labels for ``success''
and ``failure'' interchanged. To avoid this asymmetry, $u$ must
satisfy $\rho^{2}(u(y),\ell')=\rho^{2}(u(n-y),\ell')$ which $u=y$
does.

\subsection{Normal Location Family}

For a function $f:\mathcal{X}\rightarrow\mathbb{R}$ such that $f>0$
a.e. and $\int fdx=1$, the set $\left\{ f(x-\theta):\theta\in\mathbb{R}\right\} $
defines a location family indexed by $\theta$. All members of this
family have the same shape so, without loss of generality, we choose
the function $f$ having $\int_{-\infty}^{0}fdx=1/2$ so that $\theta$
is the median of $f(x-\theta$). With suitable smoothness on $f$
this set becomes our manifold $M=\left\{ f(x-\theta):\theta\in\mathbb{R}\right\} $.

For the normal location family with $\sigma>0$ known, $f=(2\pi)^{-1/2}\exp(-x^{2}/2)$
and the manifold is 
\[
M_{\sigma}=\left\{ f\left((x-\theta)/\sigma\right)/\sigma:\theta\in\mathbb{R}\right\} .
\]
For this family $\theta$ is also the mean and $y=\bar{x}=n^{-1}\sum x_{i}$
is sufficient so that $\mathcal{Y}=\mathbb{R}$. The standardized
score is $\hat{s}(\theta)=\sqrt{n}/\sigma(\bar{x}-\theta)$, a line
with slope $-\sqrt{n}/\sigma$ intersecting the horizontal axis at
$\bar{x}$. 

Figures for the normal family corresponding to those for the Bernoulli,
Figures \ref{fig:BinomPoint} and \ref{fig:BinomGen}, would have
a continuum of parallel lines in which case all vertical distributions
are identical, each being $N(0,1)$. For the normal manifold, the
horizontal distribution of the point estimator \emph{can} capture
the structure of the vertical distributions that describe the continuum
of hypotheses: the horizontal distribution for ml point estimator
at the normal model with $\theta=\theta_{\circ}$ is the affine transformation
$(\sigma/\sqrt{n})N(0,1)+\theta_{\circ}=N(\theta_{\circ},\sigma^{2}/n)$.

Outside the normal model, asymptotics provide the continuum of parallel
lines described above. Asymptotic calculations for a point estimator
$\hat{\theta}$ indicate the conditions under which, for large $n$,
the vertical distributions are normal and the corresponding generalized
estimators are linear so that confidence intervals obtained from $\hat{s}$
can be approximated by $\hat{\theta}\pm z_{\alpha}\widehat{se}$.

\subsection{Cauchy Location Family}

The manifold for the Cauchy family is $M=\left\{ f_{C}(x-\theta):\theta\in\mathbb{R}\right\} $
where $f_{C}(x)=\pi^{-1}(x^{2}+1)^{-1}$ and $\mathcal{Y}=\mathcal{X}^{n}$
since there is no sufficient statistic of lower dimension. The log
likelihood is
\begin{equation}
\ell(\theta;y)=\sum\log((x_{i}-\theta)^{2}+1)-n\log\pi\label{eq:cauchyloglik}
\end{equation}
and the standardized score estimate for $y$ is 
\begin{equation}
\hat{s}(\theta)=\left(\frac{2}{n}\right)^{1/2}\sum\frac{2(x_{i}-\theta)}{(x_{i}-\theta)^{2}+1}\label{eq:cauchysehat}
\end{equation}
since the Fisher information is $n/2$. 

\citet{Fisher1925} considers estimates based on the sample median
$\tilde{\theta}$ whose density for a sample of size $2k+1$ is 

\[
f_{{\rm med}}(z;\theta)=\frac{(2k+1)!}{(k!)^{2}\pi}\left[\frac{1}{4}-\frac{1}{\pi^{2}}\arctan^{2}(z-\theta)\right]^{k}\left[1+(z-\theta)^{2}\right]^{-1}.
\]
We consider two generalized estimators associated with the sample
median,
\begin{align*}
\bar{h}_{{\rm med}} & =(\tilde{\theta}-\theta)/\sqrt{V(\tilde{\theta})}\\
\tilde{s} & =\ell'_{{\rm med}}/\sqrt{V(\ell'_{{\rm med}})}
\end{align*}
where $\ell_{{\rm med}}=\log f_{{\rm med}}$ and $\mathcal{Y}_{{\rm med}}=\text{\ensuremath{\left\{  y:y=x_{(k+1)}\right\} } }$.
Both $h_{{\rm med}}$ and $\tilde{s}$ are functions of the data reduced
to the sample median while $\hat{s}$ is a function of the full data
$y\in\mathbb{R}^{n}$. 

Table \ref{tab:Cauchy Estimators} shows the squared slope and efficiency
($\mbox{Eff}^{\Lambda}\times$100\%) for the sample median for odd
sample sizes from $1$ to $31$. The values for $\Lambda(\tilde{s})$
and $\Lambda(\tilde{\theta})$ were obtained from \citet{Rider1960}
who reports $\Lambda(\tilde{s}$) and the variance of the median ($1/\Lambda(\tilde{\theta})$).
Since the variance of $\tilde{\theta}$ does not exist for samples
of size $n=1$ and $n=3$, $\Lambda(\tilde{\theta})$ is zero but
note $\Lambda(\tilde{s})>0$. The loss of efficiency is consistent
with the sample median not being a sufficient statistic. However,
the loss of efficiency is also due to the fact that $\tilde{\theta}$
is a point estimate; $\Lambda(\tilde{s})>\Lambda(\tilde{\theta})$
even though both are a function of the data only through the sample
median $\tilde{\theta}$. 

The adjusted sample size $n(\tilde{s})=n\mbox{Eff}^{\Lambda}(\tilde{s})$
shows the information loss in terms of a reduction in sample size.
For example, $\tilde{s}$ for a sample of size $n=15$ has smaller
slope than $\hat{s}$ for a sample of size 11 ($n(\tilde{s})=10.8$)
indicating that reducing the $15$ observations to the sample median
results in a loss of more than 4 observations. Using point estimator
$\tilde{\theta}$ rather than generalized estimator $\tilde{s}$ (both
having domain $\mathcal{Y}_{{\rm med}}$) results in the loss of an
additional observation ($n(\tilde{\theta})=9.8$). 

\begin{table}
\begin{tabular}{rrrrrrrr}
n & $\Lambda(\tilde{\theta})$ & $\Lambda(\tilde{s})$ & $\Lambda(\hat{s})$ & Eff$(\tilde{\theta})$ & Eff$(\tilde{s})$ & n$(\tilde{\theta})$ & n$(\tilde{s})$\tabularnewline
\hline 
\hline 
1 & 0 & .50000 & .5 & 0 & 100 & 0 & 1\tabularnewline
3 & 0 & 1.09064 & 1.5 & 0 & 72.71 & 0 & 2.2\tabularnewline
5 & 0.81883 & 1.74552 & 2.5 & 32.75 & 69.82 & 1.6 & 3.5\tabularnewline
7 & 1.63377 & 2.44042 & 3.5 & 46.68 & 69.73 & 3.3 & 4.9\tabularnewline
9 & 2.44703 & 3.16164 & 4.5 & 54.38 & 70.26 & 4.9 & 6.3\tabularnewline
11 & 3.25942 & 3.90109 & 5.5 & 59.26 & 70.93 & 6.5 & 7.8\tabularnewline
13 & 4.07130 & 4.65369 & 6.5 & 62.64 & 71.60 & 8.1 & 9.3\tabularnewline
15 & 4.88286 & 5.41608 & 7.5 & 65.10 & 72.21 & 9.8 & 10.8\tabularnewline
17 & 5.69418 & 6.18596 & 8.5 & 66.99 & 72.78 & 11.4 & 12.4\tabularnewline
19 & 6.50538 & 6.96171 & 9.5 & 68.48 & 73.28 & 13.0 & 13.9\tabularnewline
21 & 7.31647 & 7.74214 & 10.5 & 69.68 & 73.73 & 14.6 & 15.5\tabularnewline
23 & 8.12744 & 8.52636 & 11.5 & 70.67 & 74.14 & 16.3 & 17.1\tabularnewline
25 & 8.93839 & 9.31370 & 12.5 & 71.51 & 74.51 & 17.9 & 18.6\tabularnewline
27 & 9.74925 & 10.10363 & 13.5 & 72.22 & 74.84 & 19.5 & 20.2\tabularnewline
29 & 10.56011 & 10.89574 & 14.5 & 72.83 & 75.14 & 21.1 & 21.8\tabularnewline
31 & 11.37087 & 11.68970 & 15.5 & 73.36 & 75.42 & 22.7 & 23.4\tabularnewline
\end{tabular}

\caption{The three columns to the right of the sample size $n$ are the squared
slope for the median point estimate $\tilde{\theta}$, the score based
on the median $\tilde{s}$, and the score based on the full sample
$\hat{s}$. The next two columns give the percent $\Lambda$-efficiency
and the final two columns give the effective sample sizes. \label{tab:Cauchy Estimators}}
\end{table}

Table \ref{tab:Cauchy Estimators} does not include the ml point estimator
$\hat{\theta}$. To calculate $\Lambda(\hat{\theta})$ we need the
variance of $\hat{\theta}$ and the derivative of its expectation.
Instead of these calculations, we describe $\hat{\theta}$ using a
simulation of $100,000$ samples of size $n=15$. 

\citet{EFRON1978} consider $\hat{\theta}$ and $\hat{\theta}$ conditional
on the ancillary statistic $a=(x_{(2)}-x_{(1)},x_{(3)}-x_{(2)},\ldots,x_{(n)}-x_{(n-1)})$
where $x_{(i)}$ is the $i$th order statistic. To simplify notation
we assume the observations are ordered so that $x_{i}=x_{(i)}$. Conditioning
on $a\in\mathbb{R}^{n-1}$ results in the sample space 
\[
\mathcal{Y}_{a}=\left\{ y\in\mathcal{Y}:a=(x_{2}-x_{1},x_{3}-x_{2},\ldots,x_{n}-x_{n-1})\right\} 
\]
which is a one dimensional subspace of $\mathcal{Y}\subset\mathbb{R}^{n}$.

The point estimators $\hat{\theta}$ and $\hat{\theta}|_{a}$ and
generalized estimator $\widehat{S}$ are compared using the intervals
\begin{align*}
\hat{\theta}\pm z_{\alpha}\widehat{I}^{-1/2} & =\left\{ \theta:-z_{\alpha}<\widehat{I}^{1/2}(\theta-\hat{\theta})<z_{\alpha}\right\} \\
\hat{\theta}\pm z_{\alpha}I_{{\rm obs}}^{-1/2} & =\left\{ \theta:-z_{\alpha}<I_{{\rm obs}}^{1/2}(\theta-\hat{\theta})<z_{\alpha}\right\} \\
\widehat{S}^{-1}(-z_{\alpha}^{2}) & =\left\{ \theta:\widehat{S}(\theta)>-z_{\alpha}^{2}\right\} 
\end{align*}
where $I_{{\rm obs}}=-\ell''|_{\theta=\hat{\theta}}$ and $\widehat{I}=I|_{\theta=\hat{\theta}}$
are constants. The expressions after the equality show that the estimate-plus-se
approach are linear approximations to generalized estimators using
slope $I_{{\rm obs}}$ and $\widehat{I}$, respectively. The choice
of the latter slope is based on Fisher's argument that $V(\hat{\theta}|a)$
is $1/I_{{\rm obs}}$. This interval does not represent inference
based on $(\hat{\theta},I_{{\rm obs}})$ nor conditionally on $I_{{\rm obs}}$.
The test associated $\hat{\theta}|_{a}$ orders points in $\mathcal{Y}_{a}$.
\citet{Barndorff_Nielsen_1978} comments on this confusion regarding
the role of observed information. 

We use confidence intervals $\widehat{S}^{-1}(-z_{\alpha}^{2})$ rather
than $\hat{s}^{-1}(\pm z_{\alpha})$ since the latter do not always
exist, especially for $\alpha$ small. The score does not provide
a sensible ordering on the entire sample space $\mathcal{Y}$: for
each $\theta$, $\lim_{y\rightarrow\pm\infty}\hat{s}(\theta)(y)=0$.
A reasonable ordering would have these limits tend to infinity or
at least be bounded away from zero. \citet[page 37]{McLeish2012-ad}
describe this difficulty as ``the score function is $E$-ancillary
'at infinity' for a Cauchy model''. Their approach is to modify the
score. We find it simpler to use $\widehat{S}$. 

Coverage error for confidence intervals having nominal 5\% error are
7.5\% ($\hat{\theta}$), 6.9\% ($\hat{\theta}|a$), and 5.6\% ($\widehat{S}$)
indicating the normal distribution is a better approximation to $\widehat{S}$
than to the point estimators. \citet[pages 2-3]{McLeish2012-ad} note
the normal distribution is often a better approximation to $\hat{s}$
than to its root $\hat{\theta}$. Quantile-quantile plots for 100,000
samples show the normal distribution to be a very close approximation
to the distribution of ${\rm sign}(\hat{\theta}-\theta)\widehat{S}^{1/2}$
while the distributions of $\hat{\theta}$ and $\hat{\theta}|_{a}$
show noticeable departures from normality. This relationship to the
normal approximation holds for the other vertical and horizontal distributions:
$\hat{s}$ is closely approximated the standard normal distribution
while the median $\tilde{\theta}$ is not. 

\begin{figure}
\includegraphics{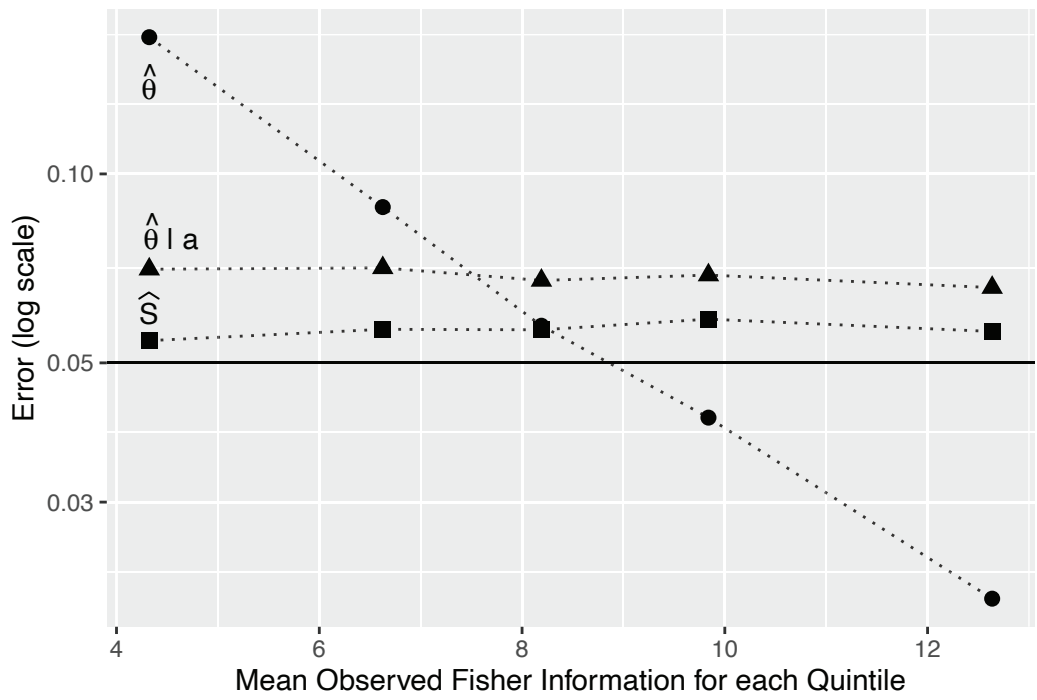}\caption{Coverage Errors and Observed Fisher Information. Simulations from
100,000 samples of size $n=15$ were ordered by the observed Fisher
information.  \label{fig:Coverage-Errors-and}}

\end{figure}

Figure \ref{fig:Coverage-Errors-and} shows that coverage error of
$\hat{\theta}$ intervals depend on $I_{{\rm obs}}$. When $I_{{\rm obs}}$
is small the coverage error is significantly larger than 5\% and when
$I_{{\rm obs}}$ is large the error significantly smaller than 5\%.
Conditioning on $a$ addresses this problem. So does $\widehat{S}$.
The generalized estimator has the advantage of avoiding other problems
that can arise with ancillary statistics: the role of approximate
ancillary statistics and choosing between multiple ancillary statistics. 

To compare the lengths we adjust the margin of errors so all intervals
have the same coverage error: $z_{.025}$ is replaced with $1.08555z_{.025}$
($\hat{\theta}$) and with $1.05518z_{.025}$ ($\hat{\theta}|_{a}$).
In keeping with our emphasis on parameter-invariance the length will
be defined using Kullback-Leibler balls (disks in this case) contained
in $M$. For $m_{\circ}\in M$, $r>0$, and $D(m_{1},m_{2})=E_{m_{1}}(\log(m_{1}/m_{2}))$
define
\[
B(m_{\circ},r)=\left\{ m\in M:D(m,m_{\circ})<r\right\} .
\]
The ball $B(m_{\circ},r)$ \emph{covers} the interval $(\theta_{lo},\theta_{hi})$
if $m\in B(m_{\circ},r)$ for all $m$ such that $\theta(m)\in(\theta_{lo},\theta_{hi})$.
The KL length of $(\theta_{lo},\theta_{hi})$ is the radius of the
smallest KL ball that covers $(\theta_{lo},\theta_{hi})$
\[
\mbox{Len}_{K\!L}(\theta_{lo},\theta_{hi})=\inf\left\{ r:\mbox{ for some }m_{\circ}\in M,B(m_{\circ},r)\mbox{ covers }(\theta_{lo},\theta_{hi})\right\} .
\]

For the Cauchy family the KL divergence is
\[
D(m_{1},m_{2})=\log\left(\left(\theta(m_{1})-\theta(m_{2})\right)^{2}+4\right)-\log4
\]
so that the KL length of the interval $(\theta_{lo},\theta_{hi})$
is a function of $|\theta_{hi}-\theta_{lo}|$. The location-scale
Cauchy family also has a closed form for the KL divergence (see \citet{2019arXiv190510965C}). 

The mean KL lengths are 0.141 ($\hat{\theta}$), 0.140 ($\hat{\theta}|_{a}$),
and 0.135 ($\widehat{S}$). Conditioning on $a$ provides little improvement
in terms of mean KL length. It might be surprising that the mean length
for $\hat{\theta}|_{a}$ is not closer to that of $\widehat{S}$ since
conditioning on $a$ recovers all the information as measured by Fisher
information. Note that information calculations for $\hat{\theta}|_{a}$
use the score for statistic $\hat{\theta}$ which differs from the
linear approximation to the score used for the $\hat{\theta}|_{a}$
confidence intervals.

\section{Discussion}

Point estimators can be made parameter invariant by defining estimators
as functions from $\mathcal{Y}$ to $M$ so that models are compared
directly to data without using any parameterization. Bias and variance
that would be expressed in a particular parameter are replaced with
KL bias and variance. Details appear in \citet{Wu2012} and \citet{Vos2015}.

The class $\mathcal{G}_{{\rm score}}$ is closely related to that
of \citet{Godambe1960}. One difference is that Godambe uses $g$
to obtain a point estimator while we consider the inferential properties
of $g$ as a function on $M$. \citet{McLeish2012-ad} emphasize the
inferential value of $g$ over its roots but still use $g$ to find
point estimators. 

The role of $\mathcal{G}_{{\rm score}}$ and $\mathcal{G}_{{\rm lrt}}$
is to provide finite sample justification for inference based on the
score and likelihood ratio using $\Lambda$-optimality. The purpose
is not to introduce new estimators. For example, the $\hat{s}$ confidence
interval (using the exact distribution) is the Clopper-Pearson interval.
Many papers report expected length for binomial confidence intervals
(e.g., \citet{Brown_2001}) but we have found none that use mean KL
length or other parameter invariant measures.

The squared slope $\Lambda$ applies to a wider class of estimators
and provides simpler comparisons than does bias and variance. If one
accepts that Fisher information describes important aspects of inference,
then the fact that $\Lambda$ and $I$ are both tensors while the
variance is not, strongly suggests that squared slope $\Lambda$ should
be used to compare estimators, thereby supporting Fisher's claim regarding
the role of $I$ in finite samples. 

Our justification for this claim relies on geometry, Fisher's view
on estimation, and his view on ``repeated sampling from the same
population.'' \citet{Fisher1955} said this phrase is foreign to
his point of view. According to \citet{Cox_2016}
\begin{quote}
{[}Fisher{]} strongly emphasized that when probability was used to
describe what underlay a set of data, he did not have in mind probability
as a limiting frequency over a large number of repetitions. Rather,
by probability Fisher meant a proportion in a hypothetical infinite
population, the data being regarded as a random sample from that hypothetical
population. 
\end{quote}
Each vertical distribution for a generalized estimator, represented
in Figure \ref{fig:BinomGen}, but common to all generalized estimators,
corresponds to Fisher's ``hypothetical infinite population''. 

Measure theory provides a rigorous treatment of Fisher's idea of a
proportion in an infinite population. This is in contrast to the notion
of repeated sampling from the same population which is not part of
mathematics. Rather, it is one interpretation for a measure theoretic
probability. \citet{Vos2022} argue that making repeated sampling
one of several possible interpretations rather than a defining aspect
of frequentist inference can address the confusion regarding p-values.
It seems the same is true for comparison of estimators: ``repeated
sampling from the same population'' encourages one to think of a
single distribution, $\mathfrak{m}$ perhaps, where bias and variance
are easily understood, but hinders the understanding of geometric
quantities such as $I$  and $\Lambda$ that are not defined for
an isolated distribution. The full utility, simplicity, and beauty
of geometric quantities, such as Fisher information and slope, are
revealed once the notion of 'hypothetical repeated sampling' is no
longer a defining property of frequentist inference.  

\bibliographystyle{apalike2}
\bibliography{vos}

\end{document}